\documentclass{amsart}

\usepackage{wrapfig}
\usepackage{epsfig}
\usepackage[dvipsnames]{xcolor}
\usepackage{caption, subcaption}

\usepackage{amsmath, amssymb, amsthm, stmaryrd, nicematrix}

\usepackage{tikz}
\usetikzlibrary{matrix,arrows,backgrounds,shapes.misc,shapes.geometric,patterns,calc,positioning, knots}
\usepackage{tikz-cd} 

\usepackage[margin=1.3in]{geometry}

\usepackage[color,matrix,arrow]{xy}
\xyoption{poly}
\xyoption{2cell}
\xyoption{all}

\usepackage[shortlabels]{enumitem}

\usepackage[colorlinks=true, pdfstartview=FitV, linkcolor=purple, citecolor=blue,
filecolor=violet, urlcolor=violet]{hyperref}

\newtheorem{thm}{Theorem}[section]
\newtheorem{prop}[thm]{Proposition}
\newtheorem{lemma}[thm]{Lemma}
\newtheorem{corollary}[thm]{Corollary}

\newtheorem{porism}[thm]{Porism}

\theoremstyle{definition}
\newtheorem{definition}[thm]{Definition}

\theoremstyle{remark}
\newtheorem{remark}[thm]{Remark}
\newtheorem{example}[thm]{Example}

\numberwithin{equation}{section}

\usepackage{subfiles, verbatim}

\newcommand{\cals}{\mathcal{S}}

\newcommand{\cala}{\mathcal{A}}

\newcommand{\za}{\alpha}

\newcommand{\zG}{\Gamma}

\newcommand{\kb}{\Bbbk}

\DeclareMathOperator{\Irr}{Irr}

\newcommand{\Gr}{\operatorname{Gr}}

\DeclareMathSymbol{\shortminus}{\mathbin}{AMSa}{"39}

\newcommand{\link}{
    \strand [thick] (1)
        to [out=45, in=135, "$5$"] (2)
        to [out=-45, in=180, "$4$"'] (downc)
        to [out=0, in=225, "$3$"', near start] (3)
        to [out=45, in=0, "$2$"'] (up)
        to [out=180, in=135, "$1$"'] (1)
        to [out=-45, in=225, "$10$"'] (2)
        to [out=45, in=180, "$9$"] (upc)
        to [out=0, in=135, "$8$", near start] (3)
        to [out =-45, in=0, "$7$"] (down)
        to [out=180, in=225, "$6$"] (1);
    \strand[thick] (1.75, 1)
        to [out = 135 , in = 45] (up)
        to [out = 225, in = 90, "$11$"', near end] (upc)
        to [out = -90, in = 90, "$12$"'] (downc)
        to [out = -90, in = 135, "$13$"', near start](down)
        to [out = -45, in = 225] (1.75,-1)
        to [out = 45, in = -90](right)
        to [out = 90, in = -45, "$14$"', pos=0] (1.75, 1);
    \flipcrossings{10,13,16,17,18, 19};
}

\title{Knot theory and cluster algebra III: Posets}
\author{Véronique Bazier-Matte}
\thanks{The first author was supported by the Discovery Grants program from the Natural Sciences and Engineering Research Council of Canada and the Research Support for New Academics from the Fonds de Recherche du Québec Nature et Technologie}
\address{Département de mathématiques et de statistique, Université Laval, Québec (Québec), G1V 0A6, Canada}
\email{veronique.bazier-matte.1@ulaval.ca}
\author{Ralf Schiffler}
\thanks{The second author was supported by the National Science Foundation grants DMS-2054561 and DMS-2348909}
\address{Department of Mathematics, University of Connecticut, Storrs, CT 06269-1009, USA}
\email{schiffler@math.uconn.edu}

\begin{document}

\begin{abstract}
    In previous work, we associated a module $T(i)$ to every segment $i$ of a link diagram $K$ and showed that there is a poset isomorphism between the submodules of $T(i)$ and the Kauffman states of $K$ relative to $i$.  In this paper, we show that the posets are distributive lattices and give explicit descriptions of the join irreducibles in both posets. We also prove that the subposet of join irreducible Kauffman states is isomorphic to the poset of the coefficient quiver of $T(i)$.
\end{abstract}

\maketitle

\section{Introduction}
\label{sect 1}
In his book \cite{K}, Kauffman introduced the Kauffman states of a link diagram $K$ with respect to a segment $i$ and used them to give a state sum formula for the Alexander polynomial of the link. 

In \cite{BMS,BMS2}, we use this approach to establish a connection to cluster algebras and representation theory. We associate a quiver $Q$ to the link diagram $K$ and consider its cluster algebra $\cala(Q)$. In \cite{BMS}, we define a representation $T(i)$ for every segment $i$ of $K$ and show that its $F$-polynomial $F_{T(i)}$ specializes to the Alexander polynomial. Moreover the specialization does not depend on the choice of the segment $i$. We also show that the poset of Kauffman states is isomorphic to the poset of subrepresentations of $T(i)$.

In \cite{BMS2}, we prove that the cluster algebra $\cala(Q)$ contains a cluster $\mathbf{x}$, the \emph{knot cluster}, in which every cluster variable specializes to the Alexander polynomial. Again the specialization is independent of the choice of the cluster variable in the knot cluster. In fact, we show that the $F$-polynomials of the knot cluster are precisely the $F$-polynomials of the representations $T(i)$, where $i$ runs over all segments of the diagram $K$.

In this paper, we use the relation to representation theory to study the combinatorial properties of the two posets. Our first main result is the following. 
\begin{thm}\label{thm intro 1}
The poset of submodules of $T(i)$ is a finite distributive lattice.
\end{thm}
We prove this theorem using an explicit $k$-vector space basis of $T(i)$ that already appeared in \cite{BMS}. Thanks to the aforementioned isomorphism of posets, we obtain a new proof of the following corollary.
\begin{corollary}
     \label{cor intro 1}
 The poset of Kauffman states relative to a segment $i$ is a finite distributive lattice. 
    \end{corollary}
This result was also proved in \cite{HK} using Tutte matchings on  triads in the sphere. Another proof was given very recently in \cite{MMSBV} who give a different realization of our connection to cluster algebras via perfect matchings of a bipartite graph that was associated to the knot in \cite{CDR}.

A famous theorem in combinatorics says that every finite distributive lattice $L$ is isomorphic to the lattice of order ideals of the subposet of join irreducible elements of $L$. We give an explicit description of the join irreducibles in terms of the subrepresentations of $T(i)$ in Proposition~\ref{prop:PropertiesLattice} and in terms of Kauffman states in Theorem~\ref{thm Kstate}. 
We show that each join irreducible subrepresentation is generated by a single element, and the collection of these generators forms a basis of $T(i)$ as a $\kb$-vector space. 

For our second main result, we need the notion of a coefficient quiver of a representation, which is a standard visualization tool. Its vertices are given by a fixed basis of the representation, and its arrows are determined by the action of the algebra on this basis. 
\begin{thm}\label{thm intro 2}
The poset of join irreducible Kauffman states relative to the segment $i$ is isomorphic to the poset of the coefficient quiver of the representation $T(i)$.   
\end{thm}

Distributive lattices appeared in the context of cluster algebras also in the following papers. \cite{CanakciSchroll} realized distributive lattices of perfect matchings of snake graphs inside the submodule lattice of corresponding string modules by restricting to the subposet of canonically embedded submodules. More generally, \cite{PilaudReadingSchroll} considered distributive lattices arising in cluster algebras from surfaces to give a poset formula for the cluster variables. These works have some similarity but are different from ours. 
The cluster algebras in this paper are not of surface type, the modules are not string modules and we consider the poset of \emph{all} submodules.

The paper is organized as follows. After a brief preliminary section, we prove Theorem~\ref{thm intro 1} in section~\ref{sect 3}. In section~\ref{sect 4} we describe the join irreducible elements as subrepresentations of $T(i)$ and prove Theorem~\ref{thm intro 2}. The description of the join irreducible Kauffman states is given in section~\ref{sect 5}. 


\section{preliminaries}
\label{sect 2}
\subsection{Knots and links}
 A \emph{knot} is a subset of $\mathbb{R}^3$ that is homeomorphic to a circle.
A \emph{link with $r$ components} is a subset  of $\mathbb{R}^3$  that is homeomorphic to $r$ disjoint circles.   Thus a knot is a link with one component. Links are considered up to ambient isotopy. 
A  link is said to be \emph{prime} if it is not the connected sum of two nontrivial links.

A \emph{link diagram} $K$ is a projection of the link into the plane that is injective  except for a finite number of  double points that are called \emph{crossing points}. In addition, the diagram carries the information at each crossing point which of the two strands is on top and which is below. 

 A \emph{curl} is a monogon in the diagram. We usually assume without loss of generality that our link diagrams are without curls, because one can always remove them (by a Reidemeister I move) without changing the link.

Let $K$ be a link diagram without curls. We consider $K$ as a cellular decomposition of $S^2$ and 
denote by $K_0$ the set of crossing points, $K_1$ the set of segments and by $K_2$ the set of regions of the complement of $K$. Let $n=|K_0|$ be the number of crossing points.  Then $|K_1|=2n$, because every crossing point is incident to four segments, and every segment is incident to two crossing points. Moreover $|K_2|=n+2$, by Euler's formula. 

\subsection{Kauffman states relative to a fixed segment} 
Let $i\in K_1$ be a segment and denote the two adjacent regions $R_i$ and $R_i'$.
A \emph{marker} is 
a pair $(x,R)$ of a crossing point $x$ and a region $R$ such that $x$ is incident to $R$. 
A \emph{Kauffman state} relative to $i$ is a set of markers, such that:

\begin{itemize}
    \item each crossing point is used in exactly one marker;
    \item each region except for $R_i$, $R_i'$ is used in exactly one marker.
\end{itemize}
The regions $R_i$, $R_i'$ are used in no marker.

\begin{figure}[ht] 
\centering
\begin{subfigure}{.45\textwidth}
  \centering
  
    \begin{tikzpicture}
    \draw[thick] (0,0) -- (4,0);
    \draw[thick] (1,1) -- (1,-1);
    \draw[thick] (3,1) -- (3,-1);
    
    \coordinate[label=below left:{$x$}] (x) at (1,0);
    \coordinate[label=below right:{$y$}] (y) at (3,0);

    \coordinate[label=below right:{$\bullet$}] (m1) at (1,0);
    \coordinate[label=above left:{$\bullet$}] (m2) at (3,0);

    
    \node at (2,0.5) {$R_1$};
    \node at (2,-0.5) {$R_2$};
    
    \end{tikzpicture}
  
  \caption{State $\cals$}
  \label{fig:state1}
\end{subfigure}%
\begin{subfigure}{.45\textwidth}
  \centering
  
    \begin{tikzpicture}
    \draw[thick] (0,0) -- (4,0);
    \draw[thick] (1,1) -- (1,-1);
    \draw[thick] (3,1) -- (3,-1);
    
    \coordinate[label=below left:{$x$}] (x) at (1,0);
    \coordinate[label=below right:{$y$}] (y) at (3,0);
    
    \coordinate[label=above right:{$\bullet$}] (m1) at (1,0);
    \coordinate[label=below left:{$\bullet$}] (m2) at (3,0);
    
    \node at (2,0.5) {$R_1$};
    \node at (2,-0.5) {$R_2$};
    \end{tikzpicture}
  
  \caption{State $\cals'$}
  \label{fig:state2}
\end{subfigure}
\caption{Kauffman counterclockwise transposition from $\cals$ to $\cals'$.}\label{fig::KauffmanTransposition}
\end{figure}
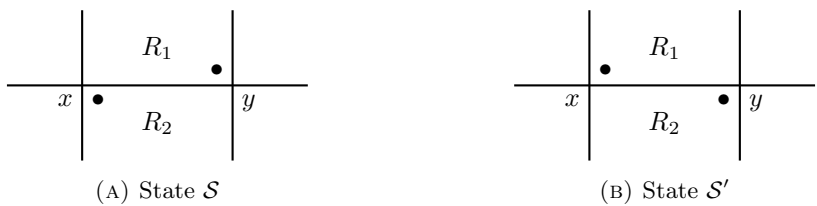

A state $\cals'$ is obtained from a state $\cals$ by a \emph{counterclockwise transposition} at a segment $j$ if $\cals'$ is obtained from $\cals$ by switching two markers at the segment $j$ as in Figure \ref{fig::KauffmanTransposition}.

More precisely, let $x$, $y$ be the endpoints of the segment $j$ and let $R_1$, $R_2$ be the adjacent regions at $j$ such that, going clockwise around $x$, we go from $R_1$ to $R_2$ crossing $j$. Then, $\cals$ contains the markers $(x, R_2)$, $(y, R_1)$, $\cals'$ contains the markers $(x, R_1)$, $(y, R_2)$ and the other markers in $\cals$ and $\cals'$ are the same.

Kauffman defines a partial order on the set of all Kauffman states by $\cals_1 < \cals_2$ if there is a sequence of counterclockwise transpositions that transforms $\cals_1$ into $\cals_2$. He proved that the resulting poset has a unique maximal and a unique minimal element. 

\subsection{Quiver representations}
Let $\kb$ be an algebraically closed field. A quiver $Q=(Q_0,Q_1)$ is an oriented graph with vertex set $Q_0$ and arrow set $Q_1$. A \emph{representation} $M=(M_i,\phi_\za)$ of $Q$ is a family of $\kb$-vector spaces $M_i$, $i\in Q_0$ and $\kb$-linear maps $\phi_\za\colon M_i\to M_j$, $\za\colon i\to j\in Q_1$. A \emph{subrepresentation} $L$ of $M$ is a family of subspaces $L_i\subset M_i$ such that $\phi_\za(L_i)\subset L_j$, for all arrows $\za\colon i\to j$. Quiver representations are equivalent to modules over the path algebra of the quiver, and we use the words representation and module interchangeably. For further details see the textbooks \cite{ASSbook,Sbook}.

\subsection{The knot module}\label{sect 2.3} We briefly recall the required material from \cite{BMS,BMS2}.  An example is given in section~\ref{sect example}. For further details, we refer to those papers. To a knot diagram $K$, 
we associate a quiver $Q=(Q_0,Q_1)$ whose vertices are the segments of $K$, thus $Q_0=K_1$. Moreover, the arrows of $Q$ are given by the following rule: \begin{itemize}
    \item [(i)] for every crossing point in $K_0$, whose segments are $a_1,a_2,a_3,a_4\in K_1$ in clockwise order, draw a clockwise oriented  4-cycle $a_1\to a_2\to a_3\to a_4\to a_1$;
    \item[(ii)] remove all 2-cycles.
\end{itemize} 

We  associate to each segment $i\in K_1$, an indecomposable representation $T(i)$ of the quiver $Q$. 
The dimension of $T(i)$ at each vertex is determined by taking successive boundaries of the knot diagram $K$ opened at the segment $i$.
Proposition 5.9 of \cite{BMS} states that if $\za\colon j\to k$ is an arrow in $Q$ then $|\dim T(i)_j -\dim T(i)_k|\le 1$.

Furthermore, there exists a basis $B(i)_j$ of $T(i)_j$, $j\in K_1$, such that  the linear map $\phi_\za$ of $T(i)$ on an arrow $\za\colon j\to k$ is of the following form, where the letter $I$ stands for an identity matrix:

\[\begin{array}{cccl}
       \begin{bNiceArray}{ccc}
        \Block{3-3}{I} \\
        & &  \\
        & &  \\
        \hline
        0 & \cdots & 0
    \end{bNiceArray}
  & or & 
\begin{bNiceArray}{c|ccc}
        0 & \Block{3-3}{ I} \\
        \vdots & & & \\
        0 & & & 
    \end{bNiceArray} 
    &\textup{if  $|\dim T(i)_j -\dim T(i)_k|= 1$;}
    \\[25pt]
    \begin{bNiceArray}{c|cc}
        0 & \Block{2-2}{I} \\
        \vdots & & \\
        \hline
        0 & \cdots & 0 
    \end{bNiceArray} 
     & or &  
  \begin{bmatrix}
      &&&& \\ &&I&&\\&&&&
  \end{bmatrix}
  & \textup{if  $\dim T(i)_j =\dim T(i)_k$.}
\end{array}\]
  In the latter case, the first matrix is a Jordan block with eigenvalue 0 and the second matrix is the identity matrix. The Jordan matrix is used if the arrow $j\to k$ corresponds to a position in the knot diagram that contains a marker of the minimal Kauffman state. Otherwise, the matrix is the identity matrix. 

The following result is our main tool in the study of the poset of Kauffman states.
\begin{thm}\cite[Theorem 6.9]{BMS}
 \label{thm lattice iso} For every segment $i$ of the link diagram $K$,
 the poset of Kauffman states of $K$ relative to the segment $i$ is isomorphic to the poset of submodules of the module $T(i)$. Moreover, two states $\cals < \cals'$ are related by a counterclockwise transposition at segment $j$ if and only if the corresponding modules satisfy $M_k'=M_k$ if $k\ne j$ and $\dim M_j'=\dim M_j+1$.
\end{thm}

\subsection{Lattices}
A  \emph{lattice} is a poset $(L,\le)$ in which each pair of elements $a,b$ admits a least upper bound $a\vee b$, called  \emph{join}, and a greatest lower bound $a\wedge b$, called \emph{meet}. A lattice $L$ is \emph{distributive} if for all $a,b,c\in L$, we have $a\wedge(b\vee c)=(a\vee b)\wedge(a\vee c).$ An element $c\in L$ is called \emph{join irreducible} if $c$ is not the minimum element and whenever $c=a\vee b$ then $c=a$ or $c=b$. In the Hasse diagram of $L$, the join irreducibles are the elements $c$ that admit a unique edge going down. This unique edge is called the \emph{descent} at $c$.

\begin{example} Figure \ref{fig:lattice} shows a lattice where the join irreducible elements are in red. \end{example}

A subset $I\subset L$ is called an \emph{order ideal} if for all $i\in I $ and $j\le i$, we have  $j\in I$. 
Given $a,b\in L$, we say $b$ \emph{covers} $a$ and write $a\lessdot b$ if $a<b$ and whenever $a<c<b$ then $c=a$ or $c=b$.

The following result is known as the Fundamental Theorem of Finite Distributive Lattices, see for example \cite[Section 3.4]{Stbook}. 

\begin{thm}
 \label{thm ftfdl}
 Any finite distributive lattice $L$ is isomorphic to the lattice of order ideals of the poset of join irreducible elements of $L$.
\end{thm}

\subsection{Coefficient quiver}
We need the following standard tool from representation theory. 
\begin{definition}
    Let $T = \left( T_j, \phi_\alpha \right)$ be a representation of a quiver $Q$ and assume $\Gamma$ has a basis $B = \bigcup_{j \in Q_0} B_j$ such that $B_j$ is a basis of $T_j$ for all vertices $j \in Q_0$. 
    For every arrow $\alpha \in Q_1$, the map $\phi_\alpha$ is given by a matrix with respect to the basis $B$.
    The \emph{coefficient quiver} $\Gamma (T, B)$ has vertex set $B$ and there is an arrow $\alpha : b \to b'$ if and only if the matrix of $\phi_\alpha$ is nonzero at position $(b,b')$.
\end{definition}

\subsection{Example}\label{sect example}
We illustrate the concepts above in an example. 
\begin{figure}[htbp]
    \centering
    \begin{subfigure}[b]{0.45\textwidth}
        \centering
        \[ \begin{tikzpicture}[scale=1.5]
            \coordinate (1) at (-1,0);
            \coordinate (2) at (0,0);
            \coordinate (3) at (1.5,0);
            \coordinate (up) at (1,1);
            \coordinate (down) at (1,-1);
            \coordinate (upc) at (0.75,0.25);
            \coordinate (downc) at (0.75,-0.25);
            \coordinate (right) at (2,0);
            \coordinate (extra) at (3,0);
            \begin{knot}[consider self intersections, ignore endpoint intersections=false, clip width=5]
            \link
            \end{knot}
        \end{tikzpicture} \]
        \caption{Link $K$}
        \label{subfig:Link}
    \end{subfigure}    
    \begin{subfigure}[b]{0.45\textwidth}
    \centering
    \begin{tikzcd}[row sep=6pt, column sep=6pt] 
        \ar[ddrrrrrr, bend left=75] \ar[d] && 11 \ar[dr] \ar[ll] && 2 \ar[dddd] \ar[ll]\\
        5 \ar[r] & 9 \ar[ru] \ar[dd] && 8 \ar[ur] \ar[dl] \\
        && 12 \ar[lu]\ar[rd] &&&& 14 \ar[lluu]\ar[lllllldd, bend left=75] \\
        10 \ar[d] & 4 \ar[l] \ar[ru] && 3 \ar[uu] \ar[ld] \\
        6\ar[uuuu, bend left = 30] \ar[rr] && 13 \ar[lu] \ar[rr] && 7 \ar[lu] \ar[rruu]
    \end{tikzcd}
    \vspace{-0.5 cm}
    \caption{Quiver $Q$}
        \label{subfig:Quiver}
    \end{subfigure}
    \begin{subfigure}[b]{0.4\textwidth}
    \centering
    \begin{tikzcd}[row sep=6pt, column sep=6pt]
        0 \ar[ddrrrrrr, bend left=75] \ar[d] && \kb \ar[dr] \ar[ll] && \kb \ar[dddd] \ar[ll]\\
        \kb \ar[r] & \kb \ar[ru] \ar[dd] && \kb^2 \ar[ur] \ar[dl] \\
        && \kb \ar[lu]\ar[rd] &&&& 0 \ar[lluu]\ar[lllllldd, bend left=75] \\
        0 \ar[d] & 0 \ar[l] \ar[ru] && 0 \ar[uu] \ar[ld] \\
        0 \ar[uuuu, bend left = 30] \ar[rr] && 0 \ar[lu] \ar[rr] && \kb \ar[lu] \ar[rruu]
    \end{tikzcd}
     \vspace{-0.5 cm}
    \caption{Representation $T(6)$}
        \label{subfig:Representation}
    \end{subfigure}
    \begin{subfigure}[b]{0.4\textwidth}
    \centering
    \begin{tikzcd}[sep=small]
        &&b_{8,2} \ar[ld] \ar[rd] \\
        b_5 \ar[rd] & b_{12} \ar[ddrr] \ar[d] && b_2 \ar[ddll] \ar[d]\\
        & b_9 \ar[d] && b_7 \ar[d] \\
        & b_{11} \ar[rd] && b_3 \ar[ld]\\
        && b_{8,1}
\end{tikzcd}
    \caption{Coefficient quiver $\Gamma (T(6),B(6))$ }
        \label{subfig:CoeffQuiver}
    \end{subfigure}
    \caption{A link, its associated quiver, the representation $T(6)$ on the Jacobian algebra given by this quiver and the coefficient quiver of this representation.}
    \label{fig:ex1}
    \end{figure}
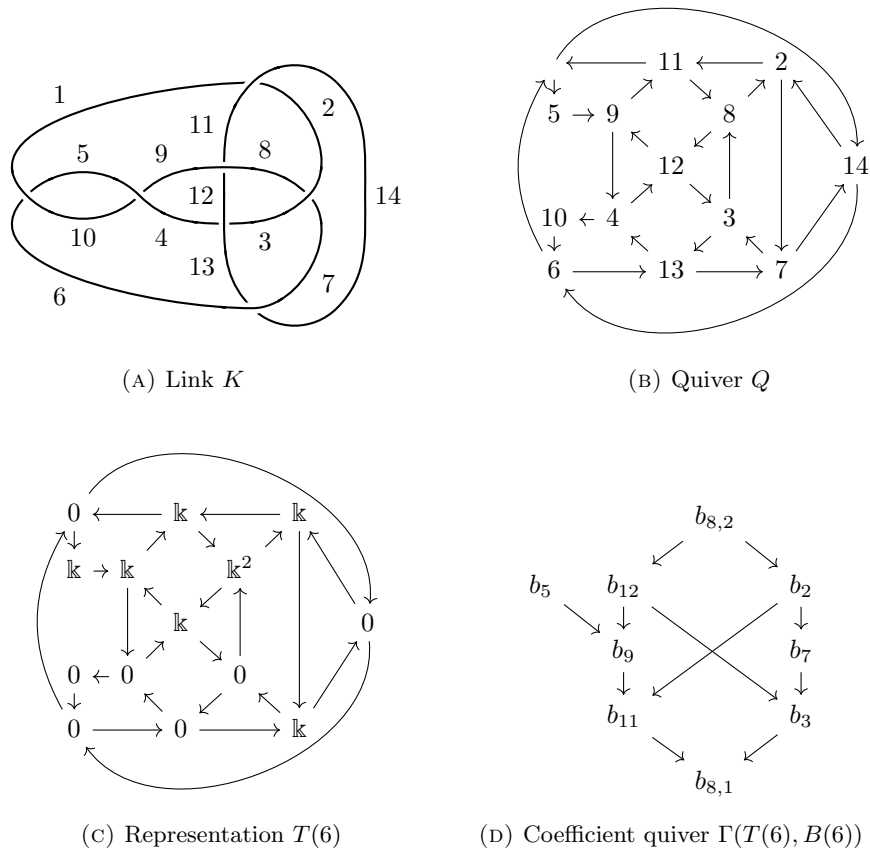

    Consider the link $K$ given in Figure \ref{subfig:Link}.
    The corresponding quiver $Q$ is shown in Figure \ref{subfig:Quiver} and the representation $T(6)$ is shown in Figure \ref{subfig:Representation},
    where the linear maps on the arrows of $T(6)$ are given by the following matrices determined by the dimensions of the vector spaces $\left[\begin{smallmatrix}
        1
    \end{smallmatrix}\right]\colon\kb\to\kb,\ \left[\begin{smallmatrix}
        0&1
    \end{smallmatrix}\right]\colon \kb^2\to\kb,\ \left[\begin{smallmatrix}
        1\\0
    \end{smallmatrix}\right]\colon \kb\to\kb^2 $.

    The coefficient quiver of $T(6)$ is given in Figure \ref{subfig:CoeffQuiver}.




\section{Distributivity}
\label{sect 3}
Let $K$ be a primitive link diagram, $i$ a segment of $K$, $T(i)$ the corresponding representation and $L(i)$ the poset of Kauffman states relative to $i$.
\begin{prop}
    The poset of Kauffman states relative to any segment $i$ is a lattice.
\end{prop}

\begin{proof}
    This follows from \cite[Theorem 6.4]{BMS}, because the poset of submodules of any module is a lattice under inclusion, where the meet is given by the intersection and the join by  the sum of two modules.
\end{proof}

\begin{example}
     Figure \ref{fig:lattice} shows the lattice of Kauffman states relative to segment $6$ in the link $K$ from Figure \ref{subfig:Link}.
    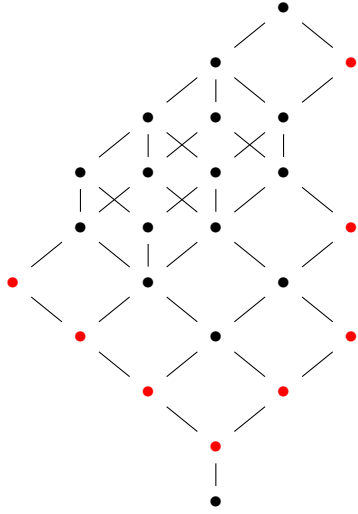
\begin{figure}[htb]
        \centering
        \[ \begin{tikzcd}[column sep = small, row sep = small]
            &&&& \bullet \ar[-, dl] \ar[-, dr]\\
            &&& \bullet \ar[-, dl] \ar[-, dr] && \textcolor{red}{\bullet} \ar[-, dl] \\
            && \bullet \ar[-, dl] \ar[-, dr] &\bullet\ar[-, dl]\ar[-, dr]\ar[-, u]& \bullet \ar[-, dl] \ar[-, d]\\
            & \bullet \ar[-, dr] \ar[-, d] &\bullet\ar[-, dl]\ar[-, dr]\ar[-, u]& \bullet \ar[-, dl] \ar[-, d] & \bullet \ar[-, dl] \ar[-, dr]\\
            & \bullet \ar[-, dl] \ar[-, dr] & \bullet \ar[-, d] & \bullet \ar[-, dl] \ar[-, dr] && \textcolor{red}{\bullet} \ar[-, dl]\\
            \textcolor{red}{\bullet} \ar[-, dr] && \bullet \ar[-, dl] \ar[-, dr] && \bullet \ar[-, dl] \ar[-, dr] \\
            & \textcolor{red}{\bullet} \ar[-, dr] && \bullet \ar[-, dl] \ar[-, dr] && \textcolor{red}{\bullet} \ar[-, dl] \\
            && \textcolor{red}{\bullet} \ar[-, dr] && \textcolor{red}{\bullet} \ar[-, dl] \\
            &&& \textcolor{red}{\bullet} \ar[-, d] \\
            &&& \bullet
        \end{tikzcd}\]
    \caption{Lattice of Kauffman states of $K$ relative to segment $6$ and of submodules of $T(6)$ for the link given in Figure \ref{subfig:Link}}
    \label{fig:lattice}
    \end{figure}
\end{example}

In general, the submodule lattice of a module is not distributive. For example, take the path algebra of the Kronecker quiver $\begin{tikzcd}[column sep=small] 1 \arrow[r, shift left] \arrow[r, shift right] & 2 \end{tikzcd}$, and the projective module at vertex $1$
\[M = P(1) = \begin{smallmatrix}
    1 \\ 2 \, 2
\end{smallmatrix}
= \begin{tikzcd}
\Bbbk \arrow[r, shift left, "{\left[\begin{smallmatrix} 1 \\0 \end{smallmatrix}\right]}"] \arrow[r, shift right, "{\left[\begin{smallmatrix} 0 \\1 \end{smallmatrix}\right]}"'] & \Bbbk^2
\end{tikzcd}.\]
Let $L$, $L'$ and $N$ three independent lines in $M_2 = \Bbbk^2$.
Then, $L, L', N$ are submodules of $M$ and $L +  L'  \supseteq N$, so $(L + L') \cap N  = N$, but $(L \cap N) + (L' \cap N) = 0 + 0 = 0$, hence the lattice is not distributive.

In the above example, the quiver Grassmannian $\Gr_{(0,1)}M$ is a projective line $\mathbb{P}^1$, which allows for the three submodules $L$, $L'$ and $N$.

Our situation is special because the module $T(i)$ has the remarkable property that, for all dimension vectors $\underline{e}$, the quiver Grassmannian $\Gr_{\underline{e}} T(i)$  is a single point or the empty set. 
Moreover, there exists a basis $B(i)$ such that the linear maps $\varphi_\alpha$, $\alpha \in Q_1$ of $T(i)$ are of the form listed in Section~\ref{sect 2}. We will show that the submodule lattice of $T(i)$ is distributive.

\begin{lemma} \label{lem:InclusionSubmod}
    Let $L$, $L'$ be submodules of $T(i)$ and denote by $d_j$, $d_j'$ their respective dimension vector at vertex $j \in Q_0$. Then, $L_j = L'_j$ if $d_j = d_j'$ and $L_j \subset L_j'$ if $d_j > d_j'$.
\end{lemma}

\begin{proof}
    According to \cite[Section 6.1]{BMS}, each vector space $(T(i))_j$ has a distinguished basis $b_{j_1}, \dots, b_{j_{d_j}}$ such that, for every chordless cycle $\omega$ from $j$ to $j$, the action of $\omega$ on the basis is given by the Jordan matrix $J_{d_j}$ with eigenvalue $0$.
    In particular, if a submodule $L$ of $T(i)$ contains the basis vector $b_{j_k}$, then it contains all the previous basis vectors $b_{j_1}, \dots, b_{j_{k-1}}$. Thus, $\dim L_j = k$ if and only if $L_j = \operatorname{span} \{b_{j_1}, \dots, b_{j_k}\}$. This completes the proof.
\end{proof}

\begin{porism} \label{por:DimSpan}
    Let $L$ be a submodule of $T(i)$.
    Then, $\dim L_j = k$ if and only if $L_j = \operatorname{span}\{b_{j_1}, \dots, b_{j_k}\}$. \qed
\end{porism}

We are now ready for our first main result. 

\begin{thm} \label{thm:DistributiveLattice}
    The submodule lattice $L(i)$ of $T(i)$ is distributive.
\end{thm}

\begin{proof}
    Let $L$, $L'$ and $L''$ be submodules of $T(i)$ and let $j \in Q_0$. Lemma \ref{lem:InclusionSubmod} implies that the vector spaces $L_j$, $L_j'$ and $L_j''$ are linearly ordered by their respective dimension $d_j$, $d_j'$ and $d_j''$. Thus, \begin{align*}
        (L_j + L'_j) \cap L_j''
         &= \min \left(\max (L_j, L_j'), L''_j\right) \\
         &= \max \left( \min (L_j, L''_j), \min (L_j, L'_j) \right) \\
         &= (L_j \cap L_j'') + (L_j' \cap L_j'') \\
         &= \left( (L \cap L'') + (L' \cap L'')\right)_j
    \end{align*}
    This shows distributivity at every vertex $j$. The statement follows because intersection and sum are additive.
\end{proof}

Combining Theorem \ref{thm:DistributiveLattice} and Theorem~\ref{thm lattice iso} yields the following result.

\begin{corollary}
    The poset of Kauffman states is a finite distributive lattice. \qed
\end{corollary}


\section{Join irreducibles}
\label{sect 4}
As before, let $K$ be a link diagram, $i$ a segment of $K$, $T(i)$ the corresponding representation and $L(i)$ the lattice of Kauffman states relative to $i$.
According to the Fundamental Theorem of Finite Distributive Lattices, the lattice $L(i)$ is isomorphic to the lattice of order ideals of the subposet of join irreducible elements of $L(i)$. In this section, we describe the join irreducibles as subrepresentations of $T(i)$.

Let $B(i) = \bigcup_{j \in Q_0} \{j_1, \dots, j_{d_j} \}$ denote the distinguished basis of $T(i)$ and let $M(j,k)$ the submodule of $T(i)$ generated by $\{b_{j_1}, \dots, b_{j_k} \}$ for $j \in Q_0$ and $k \leq d_j$.
Thus, $M(j,k)$ is generated by the first $k$ basis vectors of $T(i)_j$.

\begin{remark}\label{rem cyclic module}
It follows from the proof of Lemma~\ref{lem:InclusionSubmod} that $M(j,k)$ is generated by the single element $b_{j,k}$, because for every chordless cycle $c$ in $Q$, we have $\phi_c(b_{j,k})=b_{j,k-1}$.
\end{remark}

\begin{prop} \label{prop:PropertiesLattice}\begin{enumerate}[\rm (a)]
    \item The join irreducibles of $L(i)$ correspond to the modules $M(j,k)$ under the isomorphism between $L(i)$ and the submodule lattice of $T(i)$. \label{prop:part:a}
    \item The unique descent at $M(j,k)$ is labeled $j$. \label{prop:part:b}
    \item The join irreducibles with unique descent $j$ form a linear subposet of $L(i)$: \[M(j,1) \subset M(j,2) \subset \dots \subset M(j,d_j).\] \label{prop:part:c}
\end{enumerate}
    
\end{prop}

\begin{proof}
    By definition, $M(j,k)$ is the smallest submodule of $T(i)$ that contains $b_{j_1}, \dots, b_{j_k}$. 
    Thus, every proper submodule of $M(j,k)$ is of strictly smaller dimension at vertex $j$.
    Therefore, in the lattice, $M(j,k)$ has a descent labeled $j$.
    Moreover, $M(j,k)$ has a unique descent, because otherwise, the Hasse diagram would contain a subdiagram of the following form, with $B$ a proper submodule of $M(j,k)$:
    \[ \begin{tikzcd}[sep=small]
        & M(j,k) \arrow[-, "j"']{ld} \arrow[-, "j'"]{rd}\\
        A \arrow[-, "j'"']{rd} && B \arrow[-, "j"]{ld} \\
        &C
    \end{tikzcd}\]
    Then $B$ would have the same dimension of $M(j,k)$ at vertex $j$, which is impossible.
    Thus, $M(j,k)$ is join irreducible with unique descent labeled $j$.

    Now, let $N$ be any submodule of $T(i)$ that has dimension $k$ at vertex $j$. By Porism \ref{por:DimSpan}, $N$ contains $b_{j_1}, \dots, b_{j_k}$, and, therefore, $M(j,k)$ is a submodule of $N$.
    In particular, if $N$ is also join irreducible with unique descent labeled $j$, then $N = M(j,k)$.

    This shows \ref{prop:PropertiesLattice} \ref{prop:part:a} and \ref{prop:part:b}. Part \ref{prop:part:c} follows directly from the definition of $M(j,k)$.
\end{proof}

Denote by $\Irr(i)$ the set of join irreducibles in $L(i)$.

\begin{corollary} \label{cor:BijJoinIrrBasis}
    The map 
    $f_i\colon \Irr(i)\to B(i),\  M(j,k)\mapsto b_{j,k}$ 
    is a bijection.
\end{corollary}

\begin{proof}
    By definition, $M(j,k)$ is the submodule generated by $b_{j_1}, \dots, b_{j_k}$. Thus, $f_i$ is well-defined and bijective.
\end{proof}

Recall that the basis $B(i)$ is the vertex set of the coefficient quiver $\Gamma(T(i), B(i))$ of $T(i)$.
Our next goal is to show that the cover relations in the poset of join irreducibles are almost in bijection with the arrows in the coefficient quiver.
We start with two preliminary lemmata.

\begin{lemma} \label{lem:CoveringRlt}
    Let $M(j,k) \lessdot M(j',k')$ be a covering relation in the poset of join irreducibles of $L(i)$. Then
    \begin{enumerate}[(a)]
        \item there is an arrow $j' \to j$ in the quiver $Q$. \label{part:CoverRltA}
        \item $k = k'$ or $k+1 = k'$. More precisely, if the minimal Kauffman state has a marker between the segments 
        $j$ and $j'$, then $k+1 = k'$. Otherwise, $k = k'$. \label{part:CoverRltB}
    \end{enumerate}
\end{lemma}

\begin{proof}
    The modules $M(j,k)$, $M(j',k')$ correspond to two Kauffman states whose unique descents are given by a Kauffman transposition at the segments $j$ and $j'$ respectively.
    \begin{enumerate}[(a)]
        \item The covering relation $M(j,k) \lessdot M(j',k')$  implies that the transposition at the segment $j'$ that leads up to $M(j',k')$ must move one of the two markers involved in the transposition at the segment $j$ that leads up to $M(j,k)$.
        Since the transpositions move the markers in the counterclockwise direction around the crossing points, we deduce that $j'$ must be one of the two segments $j_1'$, $j_2'$ in Figure \ref{fig:MarkersMjk}.
        
        \begin{figure}[htb]
            \centering
            \begin{tikzpicture}[scale=0.75]
                \draw[thick] (0,0) -- (4,0)  node[midway, above] {$j$};
                \draw[thick] (1,1) -- (1,-1)  node[near start, left] {$j_1'$};
                \draw[thick] (3,1) -- (3,-1)  node[near end, right] {$j_2'$};
                \coordinate[label=above right:{$\bullet$}] (m1) at (1,0);
                \coordinate[label=below left:{$\bullet$}] (m2) at (3,0);
            \end{tikzpicture}
            \caption{Markers at $M(j,k)$} \label{fig:MarkersMjk}
        \end{figure}
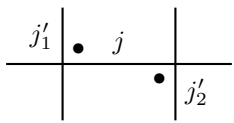
        If both $j_1'$ and $j_2'$ form a bigon with $j$ then, since $K$ is prime, $K$ must be a Hopf link and $M(j,k)=M(j,1)$ is the only join irreducible in $L(i)$. If only $j_1'$ forms a bigon with $j$, then  the transposition at $j_2'$ must be performed before the transposition at $j_1'$ is possible. Thus we may assume without loss of generality that $j'$ and $j$ do not form a bigon in $K$.
        Therefore, by definition of the quiver $Q$, there is an arrow $j' \to j$. 
        \item In the first case, each transposition at segment $j$ is preceded by a transposition at segment $j'$: 
        \[ \begin{tikzpicture}[scale=0.75]
            \draw[thick] (-1,0) -- (1,0) node[near start, above] {$j$};
            \draw[thick] (0,-1) -- (0,1) node[near start, right] {$j'$};
            \coordinate[label=below left:{$\bullet$}] (m1) at (0,0);
        \end{tikzpicture}.\]
        Thus, $\dim M(j,k)_{j'} = k$ and $\dim M(j',k')_{j'} = k+1$, and hence $k'=k+1$.

        In the second case, each transposition at segment $j'$ is preceded by a transposition at segment $j$:
        \[ \begin{tikzpicture}[scale=0.75]
            \draw[thick] (-1,0) -- (1,0) node[very near start, above] {$j$};
            \draw[thick] (0,-1) -- (0,1) node[very near start, right] {$j'$};
            \coordinate[label=below right:{$\bullet$}] (m1) at (0,0);
        \end{tikzpicture}
        \quad,\quad
        \begin{tikzpicture}[scale=0.75]
            \draw[thick] (-1,0) -- (1,0) node[very near start, above] {$j$};
            \draw[thick] (0,-1) -- (0,1) node[very near start, right] {$j'$};
            \coordinate[label=above right:{$\bullet$}] (m1) at (0,0);
        \end{tikzpicture}
        \quad,\quad
        \begin{tikzpicture}[scale=0.75]
            \draw[thick] (-1,0) -- (1,0) node[very near start, above] {$j$};
            \draw[thick] (0,-1) -- (0,1) node[very near start, right] {$j'$};
            \coordinate[label=above left:{$\bullet$}] (m1) at (0,0);
        \end{tikzpicture}.\]
        Thus, $\dim M(j,k)_{j'} = k-1$ and $\dim M(j',k')_{j'} = k$, and hence $k'=k$.\qedhere
    \end{enumerate}
\end{proof}

We are now ready for our second main result.

\begin{thm}
    The quiver of the Hasse diagram of the poset of join irreducibles of $L(i)$ is isomorphic to the subquiver of the coefficient quiver $\Gamma_{T(i)}$ obtained by removing all arrows $\alpha$ with the property that there exists a path $w \neq \alpha$ such that $s(\alpha) = s(w)$ and $t(\alpha) = t(w)$.
\end{thm}

\begin{proof}
    Let $\Gamma_{L(i)}$ be the quiver obtained from the Hasse diagram of the poset of join irreducibles of $L(i)$ by orienting each edge towards its smallest endpoint.
    We already know from Corollary \ref{cor:BijJoinIrrBasis} that there is a bijection $f_i$ between the vertices of the quivers given by mapping the Kauffman state corresponding to $M(j,k)$ to the basis element $b_{j_k}$. In order to show that this map induces a bijection on the arrows, assume we have an arrow $M(j',k') \to M(j,k)$ in $\Gamma_{L(i)}$.
    Recall that this means that $M(j,k) \lessdot M(j',k')$ in the poset of $\Irr(i)$.

    Lemma \ref{lem:CoveringRlt} \ref{part:CoverRltA} implies that there is an arrow $\alpha: j' \to j$ in  the quiver $Q$. To show that there is an arrow $b_{j',k'} \to b_{j,k}$ in $\Gamma_{T(i)}$, we must show that the matrix of $\phi_\alpha$ is nonzero at position $(b_{j',k'}, b_{j,k})$. 
    We shall use Lemma \ref{lem:CoveringRlt} \ref{part:CoverRltB} following the two distinct cases.

    Assume first the minimal Kauffman state has a marker between $j$ and $j'$ and thus $k+1=k'$. 
    We have seen in section~\ref{sect 2.3} that $\dim T(i)_j = \dim T(i)_j'$ or $\dim T(i)_j + 1 = \dim T(i)_j'$.
    In the former case, the matrix of $\phi_\alpha$ is the Jordan matrix \[\begin{bNiceArray}{c|cc}
        0 & \Block{2-2}{I} \\
        \vdots & & \\
        \hline
        0 & \cdots & 0 
    \end{bNiceArray}\]
    because the arrow 
    $\alpha$ corresponds to the position of a marker in the minimal Kauffman state.
    In the latter case, the matrix of $\phi_\alpha$ is \[\begin{bNiceArray}{c|ccc}
        0 & \Block{3-3}{I} \\
        \vdots & & & \\
        0 & & & 
    \end{bNiceArray}\]
    Figure \ref{fig:CoefQuiver} depicts the coefficient quiver locally.
    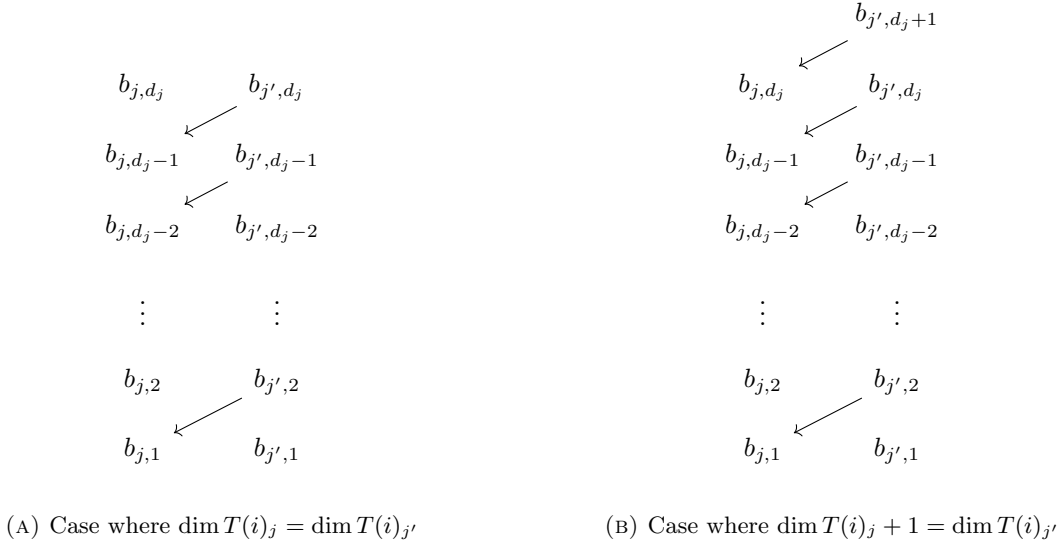
\begin{figure}[htbp]
    \centering
    \begin{subfigure}[b]{0.45\textwidth}
        \centering
        \begin{tikzcd}[sep = small]
            b_{j,d_j} & b_{j',d_j} \ar{dl} \\
            b_{j, d_j-1} & b_{j',d_j-1}\ar{dl}\\
            b_{j, d_j-2} & b_{j',d_j-2}\\
            \vdots & \vdots \\
            b_{j, 2} & b_{j',2} \ar{dl}\\
            b_{j, 1} & b_{j',1}\\
        \end{tikzcd}
        \caption{Case where $\dim T(i)_j = \dim T(i)_{j'}$}
        \label{fig:subfig1}
    \end{subfigure}
    \hfill
    \begin{subfigure}[b]{0.45\textwidth}
        \centering
        \begin{tikzcd}[sep = small]
            & b_{j',d_j+1} \ar{dl} \\
            b_{j,d_j} & b_{j',d_j} \ar{dl} \\
            b_{j, d_j-1} & b_{j',d_j-1}\ar{dl}\\
            b_{j, d_j-2} & b_{j',d_j-2}\\
            \vdots & \vdots \\
            b_{j, 2} & b_{j',2} \ar{dl}\\
            b_{j, 1} & b_{j',1}\\
        \end{tikzcd}
        \caption{Case where $\dim T(i)_j + 1 = \dim T(i)_{j'}$}
        \label{fig:subfig2}
    \end{subfigure}
    \caption{Local configuration of coefficient quiver $\Gamma_{T(i)}$ if there is a marker between segments $j$ and $j'$.}
    \label{fig:CoefQuiver}
    \end{figure}
    Thus, the only cases where there is no arrow from $b_{j',k'}$ to $b_{j,k}$ in the coefficient quiver are the following. \begin{enumerate}[i)]
        \item $\dim T(i)_j = \dim T(i)_{j'}$ and $M(j',k') = M(j',1)$: In this case, Lemma \ref{lem:CoveringRlt} implies $k = k'-1 = 0$, which is impossible.
        \item $\dim T(i)_j + 1 = \dim T(i)_{j'}$ and $M(j',k') = M(j',1)$: As in the previous case, \ref{lem:CoveringRlt} implies $k = k'-1 = 0$.
        \item $\dim T(i)_j = \dim T(i)_{j'}$ and $M(j,k) = M(j,d_j)$: Thus, Lemma \ref{lem:CoveringRlt} implies $k' = k+1 = d_j+1 = d_{j'}+1$, which is impossible.
    \end{enumerate}

    Now, assume that the minimal Kauffman state has no marker between $j$ and $j'$ and thus, $k = k'$. In this case, each transposition at segment $j'$ must be preceded by a transposition at segment $j$. Again, there are two possibilities: $\dim T(i)_j = \dim T(i)_{j'}$ or $\dim T(i)_j = \dim T(i)_{j'} + 1$.
    In the former case, the matrix of $\phi_\alpha$ is the identity matrix, since the marker of the minimal Kauffman state does not have a marker in the region corresponding to the arrow $\alpha$. In particular, the coefficient quiver has an arrow $b_{j',k} \to b_{j,k}$ for all $k$,  and there is nothing to show. 
    In the latter case, the matrix of $\phi_\alpha$ is \[\begin{bNiceArray}{ccc}
        \Block{3-3}{I} \\
        & &  \\
        & &  \\
        \hline
        0 & \cdots & 0
    \end{bNiceArray}\]
    and, hence, the coefficient quiver has an arrow $b_{j',k} \to b_{j,k}$ for all $k = 1, \dots, d_{j'}$. This situation is illustrated in Figure \ref{fig:CoefQuiver2}.
    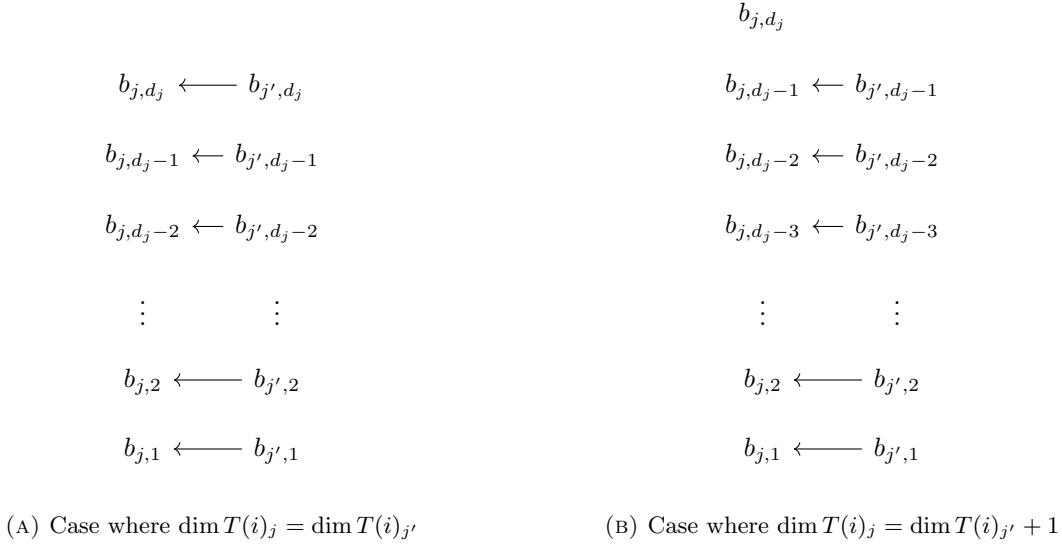
\begin{figure}[htbp]
    \centering
    \begin{subfigure}[b]{0.45\textwidth}
        \centering
        \begin{tikzcd}[sep = small]
            b_{j,d_j} & b_{j',d_j} \ar{l} \\
            b_{j, d_j-1} & b_{j',d_j-1}\ar{l}\\
            b_{j, d_j-2} & b_{j',d_j-2}\ar{l}\\
            \vdots & \vdots \\
            b_{j, 2} & b_{j',2} \ar{l}\\
            b_{j, 1} & b_{j',1} \ar{l}\\
        \end{tikzcd}
        \caption{Case where $\dim T(i)_j = \dim T(i)_{j'}$}
        \label{fig:subfig3}
    \end{subfigure}
    \hfill
    \begin{subfigure}[b]{0.45\textwidth}
        \centering
        \begin{tikzcd}[sep = small]
            b_{j,d_j} & \\
            b_{j,d_j-1} & b_{j',d_j-1} \ar{l} \\
            b_{j, d_j-2} & b_{j',d_j-2}\ar{l}\\
            b_{j, d_j-3} & b_{j',d_j-3} \ar{l}\\
            \vdots & \vdots \\
            b_{j, 2} & b_{j',2} \ar{l} \\
            b_{j, 1} & b_{j',1} \ar{l}\\
        \end{tikzcd}
        \caption{Case where $\dim T(i)_j = \dim T(i)_{j'} + 1$}
        \label{fig:subfig4}
    \end{subfigure}
    \caption{Local configuration of coefficient quiver $\Gamma_{T(i)}$ if there is no marker between segments $j$ and $j'$.}
    \label{fig:CoefQuiver2}
    \end{figure}
    There is no arrow ending at $b_{j, d_j} = b_{j, d_{j'+1}}$, but in this case, $M(j,k) = M(j,d_j)$ and Lemma \ref{lem:CoveringRlt} \ref{part:CoverRltB} implies $k' = k = d_j = d_{j'+1}$, which is impossible.
    
    This show that if there is an arrow $M(j,k) \to M(j',k')$ in the Hasse quiver $\Gamma_{L(i)}$ of the poset of join irreducibles of $L(i)$ then there is an arrow $b_{j,k} \to b_{j'k}$ in the coefficient quiver $\Gamma_{T(i)}$ of $T(i)$.

    Conversely, suppose $\gamma : b_{j',k'} \to b_{j,k}$ is an arrow in $\Gamma_{T(i)}$.
    The module $M(j',k')$ is the submodule of $T(i)$ generated by $b_{j',k'}$. The existence of the arrow $\gamma$ implies that there exists an arrow $\alpha : j' \to j$ in $Q_1$ such that the action of $\alpha$ on $T(i)$ sends the basis element $b_{j',k'}$ to the basis element $b_{j,k}$.
    Thus, $b_{j,k} \in M(j',k)$, and consequently $M(j,k)$ is a submodule of $M(j',k')$. It follows that $M(j,k) < M(j',k')$ in the poset of join irreducibles. Hence each arrow in $\Gamma_{T(i)}$ corresponds to a single arrow or a path of length at least 2 in the Hasse quiver of $\Irr(i)$.

    If there is a path $M(j',k') = M_0 \to M_1 \to \dots \to M_t = M(j,k)$ with $t\ge 2$ in the Hasse quiver, it gives arise to a path $b_{j',k'} = b_0 \to b_1 \to \dots \to b_t = b_{j,k}$ in $\zG_{T(i)}$. By the first part of the proof, the arrow $\gamma$ is parallel to this path.
    
    Otherwise,  there is an arrow $\alpha'$ in the Hasse quiver and $M(j,k) \lessdot M(j',k')$. 
    Suppose that there is as well a path $w$ parallel to $\gamma$ in the coefficient quiver. Each arrow in this path corresponds to a single arrow or a path in the Hasse diagram in such way that the whole path in $\Gamma_{T(i)}$ corresponds to a path $w'$ in the Hasse diagramm. This is a contradiction, because the arrow $\alpha'$ would be parallel to the path $w$ in the Hasse diagramm.
    This completes the proof.
\end{proof}

\begin{example}
Consider the link depicted in Figure \ref{fig:LinkNonIso}. The lattice of join irreducible Kauffman states relative to segment $9$ (Figure \ref{fig:LatticeLinkNonIso}) and the coefficient quiver of the representation $T(9)$ with basis $B(9)$ (Figure \ref{fig:CoeffQuiverLinkNonIso}) are not isomorphic.
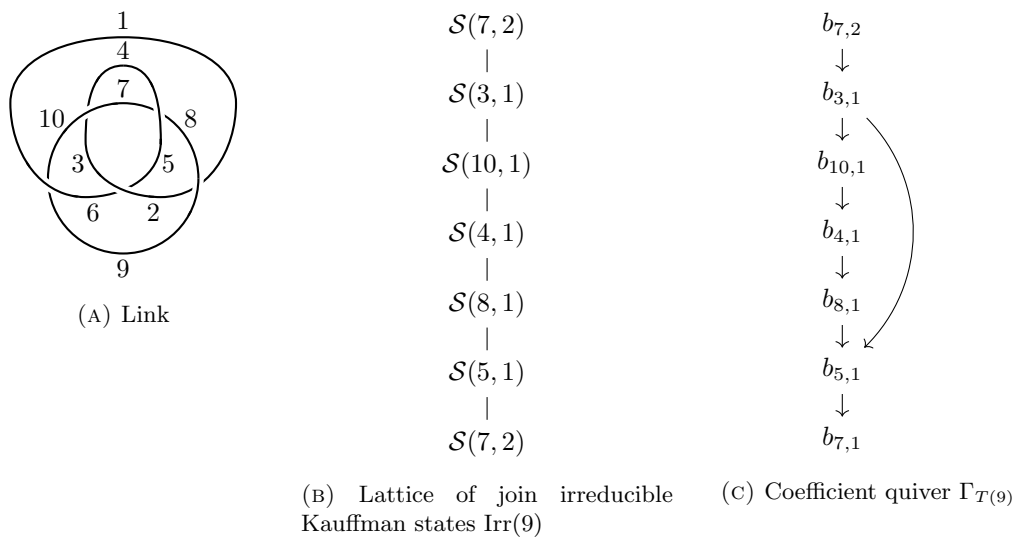
\begin{figure}[htbp]
    \centering
    \begin{subfigure}[T]{0.3\textwidth}
    \centering
    \begin{tikzpicture}
    \begin{knot}[consider self intersections, ignore endpoint intersections=false, clip width=5]
        \strand[thick] (0,0) circle[radius=1];
        \strand[thick] (0,1.5)
            to [out = 180, in = 90] (-0.5,0.5)
            to [out=270, in=180] (0.5,-0.25)
            to [out=0, in=270] (1.5,1)
            to [out=90, in=90] (-1.5,1)
            to [out=270, in=180] (-0.5,-0.25)
            to [out=0, in=270] (0.5,0.5)
            to [out=90, in=0] (0,1.5);
        \flipcrossings{1,2,4};
        \node at (0,2.1) {$1$};
        \node at (0.4,-0.45) {$2$};
        \node at (-0.6,0.2) {$3$};
        \node at (0,1.7) {$4$};
        \node at (0.6,0.2) {$5$};
        \node at (-0.4,-0.45) {$6$};
        \node at (0,1.2) {$7$};
        \node at (0.9,0.8) {$8$};
        \node at (0,-1.2) {$9$};
        \node at (-0.95,0.8) {$10$};
     \end{knot}
     \end{tikzpicture}
     \caption{Link}
     \label{fig:LinkNonIso}
     \end{subfigure}
     \begin{subfigure}[T]{0.33\textwidth}
     \centering
         \begin{tikzcd}[sep = small]
             \cals(7,2) \ar[-,d] \\
             \cals(3,1) \ar[-,d] \\
             \cals(10,1) \ar[-,d] \\
             \cals(4,1) \ar[-,d] \\
             \cals(8,1) \ar[-,d] \\
             \cals(5,1) \ar[-,d] \\
             \cals(7,2)
        \end{tikzcd}    
        \caption{Lattice of join irreducible Kauffman states $\Irr(9)$}
        \label{fig:LatticeLinkNonIso}
    \end{subfigure}
     \begin{subfigure}[T]{0.33\textwidth}
     \centering
         \begin{tikzcd}[sep = small]
             b_{7,2} \ar[d] \\
             b_{3,1} \ar[d] \ar[dddd, bend left=45] \\
             b_{10,1} \ar[d] \\
             b_{4,1} \ar[d] \\
             b_{8,1} \ar[d] \\
             b_{5,1} \ar[d] \\
             b_{7,1}
        \end{tikzcd}    
        \caption{Coefficient quiver $\Gamma_{T(9)}$}
        \label{fig:CoeffQuiverLinkNonIso}
    \end{subfigure}
    \caption{Link where the lattice of join irreducible Kauffman states and the coefficient quiver are not isomorphic}
    \label{fig:NonIsomorphic}
    \end{figure}
\end{example}

\begin{corollary}
    \label{cor irr}
    The poset of join irreducibles $\Irr(i)$ is isomorphic to the poset of the coefficient quiver $\zG_{T(i)}.$ 
\end{corollary}
\begin{proof}
    This follows from the theorem, because an edge parallel to a path is redundant  in the Hasse diagram, by transitivity.
\end{proof}


\section{The join irreducible Kauffman states}
\label{sect 5}
In this section, we describe the Kauffman states relative to segment $i$ that are join irreducible. In the previous section, we have seen that the join irreducibles correspond bijectively to the subrepresentations $M(j,k)$ of $T(i)$ that are generated by the basis element $b_{j,k}$ in $T(i)$. Here $j$ is a segment in the link diagram $K$ and $k$ is an integer between $1$ and the dimension $d_{k}$ of $T(i)$ at vertex $j$.  Thus for each segment $j$, and each $k=1,2,\ldots, d_k$, there exists a distinct Kauffman state $\cals(j,k)$ that is join irreducible in the lattice of Kauffman states $L(i)$, whose unique descent is given by the transposition at $j$ and such that every path in $L(i)$ from $\cals(j,k)$ to the minimal state uses the transposition at $j$ exactly $k$ times.

\begin{figure}
\begin{center}
\begingroup%
  \makeatletter%
  \providecommand\color[2][]{%
    \errmessage{(Inkscape) Color is used for the text in Inkscape, but the package 'color.sty' is not loaded}%
    \renewcommand\color[2][]{}%
  }%
  \providecommand\transparent[1]{%
    \errmessage{(Inkscape) Transparency is used (non-zero) for the text in Inkscape, but the package 'transparent.sty' is not loaded}%
    \renewcommand\transparent[1]{}%
  }%
  \providecommand\rotatebox[2]{#2}%
  \newcommand*\fsize{\dimexpr\f@size pt\relax}%
  \newcommand*\lineheight[1]{\fontsize{\fsize}{#1\fsize}\selectfont}%
  \ifx\svgwidth\undefined%
    \setlength{\unitlength}{74.46087502bp}%
    \ifx\svgscale\undefined%
      \relax%
    \else%
      \setlength{\unitlength}{\unitlength * \real{\svgscale}}%
    \fi%
  \else%
    \setlength{\unitlength}{\svgwidth}%
  \fi%
  \global\let\svgwidth\undefined%
  \global\let\svgscale\undefined%
  \makeatother%
  \begin{picture}(1,0.80579228)%
    \lineheight{1}%
    \setlength\tabcolsep{0pt}%
    \put(0,0){\includegraphics[width=\unitlength,page=1]{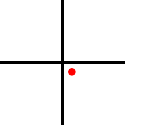}}%
    \put(0.06043441,0.44318579){\makebox(0,0)[lt]{\lineheight{1.25}\smash{\begin{tabular}[t]{l}$s$\end{tabular}}}}%
    \put(0.43311335,0.69499584){\makebox(0,0)[lt]{\lineheight{1.25}\smash{\begin{tabular}[t]{l}$t_1$\end{tabular}}}}%
    \put(0.61441664,0.28202732){\makebox(0,0)[lt]{\lineheight{1.25}\smash{\begin{tabular}[t]{l}$t_2$\end{tabular}}}}%
  \end{picture}%
\endgroup%

\caption{Four segments at a crossing point in $K$. The red bullet indicates the marker of the minimal Kauffman state.  If $s\in S$ and $t_1,t_2\in T$, then $t_1,t_2$ lie between $s$ and the minimal state marker in clockwise direction, hence they are contained in the closure of $S$ in $T$.}
\label{fig closure}
\end{center}
\end{figure}
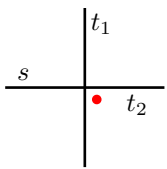

We need the following notion of closure. An illustration is given in Figure~\ref{fig closure}.
\begin{definition}
 Let $i\in K_1$ be a segment 
 and $S\subset T\subset K_1$ be 	 subsets of segments. \begin{enumerate}[(a)]
    \item  We say that $S$ is \emph{closed in $T$ under successors} (relative to $i$)  if the following condition holds.  For all segments $s\in S$, $t\in T$ that have a common  endpoint $x$, such that $t$ lies between $s$ and the minimal state marker at $x$ in clockwise direction, the segment $t$ is contained in $S$.
     \item The \emph{closure} of $S$ in $T$ is the smallest subset of $T$ containing $S$ that is closed in $T$ under successors.
\end{enumerate}
\end{definition}

\begin{definition}
 (a) Let $M(j,k)$ be a join irreducible in $K(i)$. We define a partition \[K_1=\bigsqcup_{d=0,1,\ldots,k} K(d)\] of the set of segments in $K$ as follows.
\begin{itemize}
\item 
For $d=k$ we let $K(d)$ be the closure of the set $\{j\}$ in $K_1$. 

\item Recursively, for $d=k-1, k-2,\ldots 1$, 
we let $K(d)$ be the closure in  $K_1\setminus\big(\cup_{\ell>d} K(\ell)\big)$ of the set
\[\left\{e\in K_1\setminus\big(\cup_{\ell>d} K(\ell)\big)\mid\textup{$e$ is incident to a segment of $K(d+1)$} 
\right\}\]

\item Finally, for $d=0$, we let $K(0)=K_1\setminus\big(\cup_{\ell>0} K(\ell)\big)$.

\end{itemize}

(b) If a segment $a$ lies in the set $K(\ell)$, we say $a$ is at \emph{level} $\ell$. 
\end{definition}

\begin{figure}
\begin{center}
\scalebox{0.8}{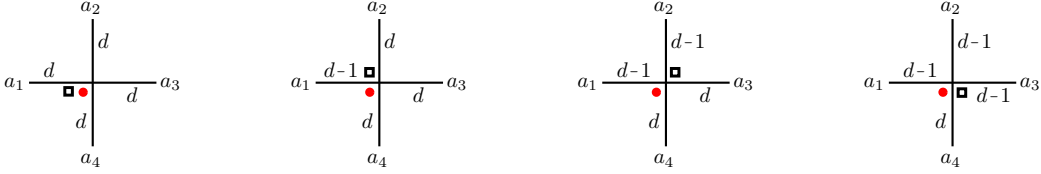}
\caption{The four cases for the markers in $\cals(j,k)$ and for the levels at a crossing point. The marker of $\cals(j,k)$ is illustrated by a black square and the marker of the minimal state by a red bullet point. The segments are labeled $a_i$ and the levels $d$ or $d-1$.}
\label{fig markers}
\end{center}
\end{figure}

We are now ready to define the Kauffman state associated to $M(j,k)$. See Figure~\ref{fig markers} for an illustration.
\begin{definition}
 Define a set of markers $\cals(j,k)$ by placing the marker at the crossing point $x\in K_0$ in the region where the level increases in clockwise direction, if such a region exists. If no such region exists, then all four segments at $x$ are at the same level, and we place the marker at the position of the minimal state marker at $x$.  
\end{definition}

\begin{example}
Let's compute $\cals(12,1)$ and $\cals(8,2)$ from the link $K$ given at Section \ref{sect example}. First, we need to know the minimal Kauffman state, which is given in Figure \ref{subfig:MinState}.
To determine the closure of $\{ 12 \}$, we consider the position of the minimal marker at its endpoints. At the top endpoint of segment 12, the segments 9, 11, and 8 lie between the segment 12 and the marker, and therefore they are contained in the closure. At the bottom endpoint of segment 12, the segment 3 lies between 12 and the marker, and thus it is contained in the closure. Next we consider the endpoints of the newly added segments. At the left endpoint of 9, the marker is directly clockwise of 9, so no segment is added here. The same is true at the endpoints of the segments 11 and 8. At the right endpoint of segment 3, we cross the segment 8 to get to the marker, and thus segment 8 must be in the closure. However, we have already added segment 8 in the previous step, and so we are done. 
Thus the closure of $\{ 12 \}$ in $K_1$ is $\{3,8,9,11,12\}$. Therefore, $K_{12,1}(1) = \{ 3,8,9,11,12 \}$ and $K_{12,1}(0) = K_1 \setminus K_{12,1}(1)$, as illustrated at Figure \ref{subfig:S12levels}, which gives in turn the Kauffman state from Figure \ref{subfig:S12}.

Now let's compute $\cals(8,1)$ and $\cals(8,2)$. The set $\{8\}$ is already closed in $K_1$, so $K(2)=\{8\}$ and the state $\cals(8,2)$ is obtained from the minimal state by the transposition at 8. 
The set $K(1)$ is the closure of the segments incident to 8 inside $K_1\setminus\{8\}$. The segments incident to 8 are $2,3,7,9,11,12$, and this set is already closed in  $K_1\setminus\{8\}$. Therefore $K(1)=\{2,3,7,9,11,12\}$, as illustrated in Figure \ref{subfig:s8levels} and the corresponding state $\cals(8,1)$ is shown in Figure \ref{subfig:S8}.
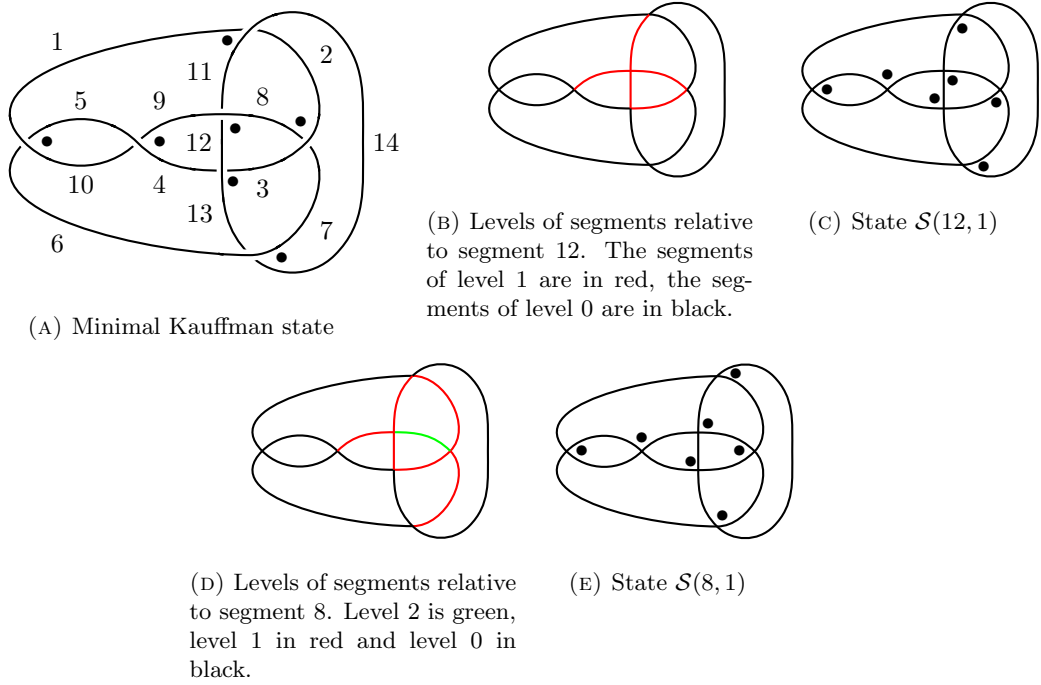
\begin{figure}
    \begin{subfigure}[t]{0.42\textwidth}
        \centering
        \[ \begin{tikzpicture}[scale=1.5]
            \coordinate (1) at (-1,0);
            \coordinate (2) at (0,0);
            \coordinate (3) at (1.5,0);
            \coordinate (up) at (1,1);
            \coordinate (down) at (1,-1);
            \coordinate (upc) at (0.75,0.25);
            \coordinate (downc) at (0.75,-0.25);
            \coordinate (right) at (2,0);
            \coordinate (extra) at (3,0);
            \begin{knot}[consider self intersections, ignore endpoint intersections=false, clip width=5]
            \link
            \end{knot}
            \node at (-0.8,0) {$\bullet$};
            \node at (0.2,0) {$\bullet$};
            \node at (1.45,0.18) {$\bullet$};
            \node at (0.8,0.9) {$\bullet$};
            \node at (1.28,-1.03) {$\bullet$};
            \node at (0.87,0.12) {$\bullet$};
            \node at (0.85,-0.35) {$\bullet$};
        \end{tikzpicture} \]
        \caption{Minimal Kauffman state}
        \label{subfig:MinState}
    \end{subfigure}
    \begin{subfigure}[t]{0.29\textwidth}
        \centering
        \[ \begin{tikzpicture}[scale=1]
            \coordinate (1) at (-1,0);
            \coordinate (2) at (0,0);
            \coordinate (3) at (1.5,0);
            \coordinate (up) at (1,1);
            \coordinate (down) at (1,-1);
            \coordinate (upc) at (0.75,0.25);
            \coordinate (downc) at (0.75,-0.25);
            \coordinate (right) at (2,0);
            \coordinate (extra) at (3,0);
            \begin{knot}[consider self intersections, ignore endpoint intersections=true]
            \strand [thick] (1)
                to [out=45, in=135] (2)
                to [out=-45, in=180] (downc);
            \strand[thick, red] (downc) 
                to [out=0, in=225] (3);
            \strand [thick] (3)
                to [out=45, in=0] (up)
                to [out=180, in=135] (1)
                to [out=-45, in=225] (2);
            \strand [thick, red] (2)
                to [out=45, in=180] (upc)
                to [out=0, in=135] (3);
            \strand [thick] (3)
                to [out =-45, in=0] (down)
                to [out=180, in=225] (1);
            \strand[thick, red] (up)
                to [out = 225, in = 90] (upc)
                to [out = -90, in = 90] (downc);
            \strand [thick] (downc)
                to [out = -90, in = 135=](down)
                to [out = -45, in = 225] (1.75,-1)
                to [out = 45, in = -90](right)
                to [out = 90, in = -45] (1.75, 1)
                to [out = 135 , in = 45] (up);
            \end{knot}
        \end{tikzpicture} \]
        \caption{Levels of segments relative to segment 12. The segments of level 1 are in red, the segments of level 0 are in black.}
        \label{subfig:S12levels}
    \end{subfigure}
    \begin{subfigure}[t]{0.25\textwidth}
        \centering
        \[ \begin{tikzpicture}[scale=1]
            \coordinate (1) at (-1,0);
            \coordinate (2) at (0,0);
            \coordinate (3) at (1.5,0);
            \coordinate (up) at (1,1);
            \coordinate (down) at (1,-1);
            \coordinate (upc) at (0.75,0.25);
            \coordinate (downc) at (0.75,-0.25);
            \coordinate (right) at (2,0);
            \coordinate (extra) at (3,0);
            \begin{knot}[consider self intersections, ignore endpoint intersections=true]
            \strand[thick] (1)
                to [out=45, in=135] (2)
                to [out=-45, in=180] (downc) 
                to [out=0, in=225] (3)
                to [out=45, in=0] (up)
                to [out=180, in=135] (1)
                to [out=-45, in=225] (2)
                to [out=45, in=180] (upc)
                to [out=0, in=135] (3)
                to [out =-45, in=0] (down)
                to [out=180, in=225] (1);
            \strand[thick] (up)
                to [out = 225, in = 90] (upc)
                to [out = -90, in = 90] (downc)
                to [out = -90, in = 135=](down)
                to [out = -45, in = 225] (1.75,-1)
                to [out = 45, in = -90](right)
                to [out = 90, in = -45] (1.75, 1)
                to [out = 135 , in = 45] (up);
            \end{knot}
            \node at (-0.8,0) {$\bullet$};
            \node at (0,0.2) {$\bullet$};
            \node at (1.45,-0.18) {$\bullet$};
            \node at (1,0.8) {$\bullet$};
            \node at (1.28,-1.03) {$\bullet$};
            \node at (0.87,0.12) {$\bullet$};
            \node at (0.63,-0.12) {$\bullet$};
        \end{tikzpicture} \]
        \caption{State $\cals(12,1)$}
        \label{subfig:S12}
    \end{subfigure}
    \begin{subfigure}[t]{0.29\textwidth}
        \centering
        \[ \begin{tikzpicture}[scale=1]
            \coordinate (1) at (-1,0);
            \coordinate (2) at (0,0);
            \coordinate (3) at (1.5,0);
            \coordinate (up) at (1,1);
            \coordinate (down) at (1,-1);
            \coordinate (upc) at (0.75,0.25);
            \coordinate (downc) at (0.75,-0.25);
            \coordinate (right) at (2,0);
            \coordinate (extra) at (3,0);
            \begin{knot}[consider self intersections, ignore endpoint intersections=true]
            \strand [thick] (1)
                to [out=45, in=135] (2)
                to [out=-45, in=180] (downc);
            \strand[thick, red] (downc) 
                to [out=0, in=225] (3);
            \strand [thick,red] (3)
                to [out=45, in=0] (up);
                \strand[thick](up)
                to [out=180, in=135] (1)
                to [out=-45, in=225] (2);
            \strand [thick, green] (upc)
                to [out=45, in=180] (upc)
                to [out=0, in=135] (3);
            \strand [thick, red] (2)
                to [out=45, in=180] (upc)
                to [out=0, in=135] (upc);
            \strand [thick, red] (3)
                to [out =-45, in=0] (down);
                \strand[thick] (down)
                to [out=180, in=225] (1);
            \strand[thick, red] (up)
                to [out = 225, in = 90] (upc)
                to [out = -90, in = 90] (downc);
            \strand [thick] (downc)
                to [out = -90, in = 135=](down)
                to [out = -45, in = 225] (1.75,-1)
                to [out = 45, in = -90](right)
                to [out = 90, in = -45] (1.75, 1)
                to [out = 135 , in = 45] (up);
            \end{knot}
        \end{tikzpicture}\] 
        \caption{Levels of segments relative to segment 8. Level 2 is green, level 1 in red and level 0 in black.}
        \label{subfig:s8levels}
    \end{subfigure}\begin{subfigure}[t]{0.25\textwidth}
        \centering
  \[ \begin{tikzpicture}[scale=1]
            \coordinate (1) at (-1,0);
            \coordinate (2) at (0,0);
            \coordinate (3) at (1.5,0);
            \coordinate (up) at (1,1);
            \coordinate (down) at (1,-1);
            \coordinate (upc) at (0.75,0.25);
            \coordinate (downc) at (0.75,-0.25);
            \coordinate (right) at (2,0);
            \coordinate (extra) at (3,0);
            \begin{knot}[consider self intersections, ignore endpoint intersections=true]
            \strand[thick] (1)
                to [out=45, in=135] (2)
                to [out=-45, in=180] (downc) 
                to [out=0, in=225] (3)
                to [out=45, in=0] (up)
                to [out=180, in=135] (1)
                to [out=-45, in=225] (2)
                to [out=45, in=180] (upc)
                to [out=0, in=135] (3)
                to [out =-45, in=0] (down)
                to [out=180, in=225] (1);
            \strand[thick] (up)
                to [out = 225, in = 90] (upc)
                to [out = -90, in = 90] (downc)
                to [out = -90, in = 135=](down)
                to [out = -45, in = 225] (1.75,-1)
                to [out = 45, in = -90](right)
                to [out = 90, in = -45] (1.75, 1)
                to [out = 135 , in = 45] (up);
            \end{knot}
            \node at (-0.8,0) {$\bullet$};
            \node at (0,0.18) {$\bullet$};
            \node at (1.30,0) {$\bullet$};
            \node at (1.25,1.03) {$\bullet$};
            \node at (1.07,-0.86) {$\bullet$};
            \node at (0.88,0.36) {$\bullet$};
            \node at (0.65,-0.15) {$\bullet$};
        \end{tikzpicture} \]        \caption{State $\cals(8,1)$}
        \label{subfig:S8}
    \end{subfigure}
    \caption{Computation of join irreducible Kauffman states}
    \label{fig:enter-label}
\end{figure}
\end{example}

The remainder of this section is devoted to the proof that $\cals(j,k)$ is join irreducible in the lattice $L(i)$ of Kauffman states relative to $i$. 

\begin{lemma}
 \label{lem 41}
 Let $x$ be a crossing point with segments $a_1,a_2,a_3,a_4$ in clockwise order such that the minimal state marker at $x$ lies between $a_4$ and $a_1$. Then the vector of the levels of $ a_1,a_2,a_3,a_4$ is of the form $(d,d,d,d),(d-1,d,d,d),(d-1,d-1,d,d), $ or $(d-1,d-1,d-1,d)$, for some $d>0$, see Figure~\ref{fig markers}. 
\end{lemma}
\begin{proof}
 This follows directly from the definition of the partition of $K_1$. 
\end{proof}

We say the level is \emph{uniform at a crossing point $x$} if all four segments at $x$ are at the same level.
\begin{lemma}
 \label{lem 42}
 Let $x$ be a crossing point of $K$. 
 
 (a) If the level is uniform at $x$, then the minimal state marker is equal to the marker in $\cals(j,k)$. 

(b) If the level is not uniform at $x$, then the minimal state marker lies in the region where the level decreases in clockwise direction, and the marker of $\cals(j,k)$ lies in the region where the level increases, see Figure~\ref{fig markers}. 
\end{lemma}

\begin{proof}
 This follows directly from the definition of the partition $\cals(j,k)$.
\end{proof}

\begin{corollary}
 \label{cor 43}
 Every region of $K$ except the two regions incident to the segment $i$ contains a unique marker of $\cals(j,k)$.
\end{corollary}
\begin{proof}
 Assume the region $R$ has two markers of $\cals(j,k)$, and let $x,y$ be the corresponding crossing points. If the level  is not uniform at $x$ and at $y$ then the level increases at both $x$ and $y$ when going around the region $R$ in a counterclockwise direction. Therefore there must be two crossing points $x',y'$ at $R$ where the level decreases. By Lemma \ref{lem 42}, this would imply that the region $R$ contains two markers of the minimal Kauffman state, which is a contradiction.
 
 Now suppose the level is uniform at $x$. Then the minimal state marker at $x$ is equal to the marker of $\cals(j,k)$, and thus it lies in the region $R$ as well.  Since $R$ cannot have two markers of the minimal state, it follows that the level is not uniform at $y$. Thus it increases at $y$. Consequently, there must be another point $y'$ at $R$ where the level decreases again, and then the minimal state marker at $y'$ lies in the region $R$. This again is a contradiction, since the minimal state would have two markers in $R$.
 
 This shows that $R$ contains at most one marker of $\cals(j,k)$. 
 Now suppose $R$ does not contain a marker of $\cals(j,k)$. Then the level is uniform in $R$. Therefore, if $R$ contains a marker of the minimal state, then by Lemma~\ref{lem 42} this marker would also be a marker of $\cals(j,k)$, a contradiction. Thus $R$ does not contain a marker of the minimal Kauffman state, and hence $R$ is a region incident to the segment $i$. 
\end{proof}

We are now ready for the main result of this section. 
\begin{thm}
 \label{thm Kstate}
 The join irreducible Kauffman states are precisely the states $\cals(j,k)$. Moreover, the isomorphism of Theorem~\ref{thm lattice iso} maps
 $\cals(j,k)$  to the representation $M(j,k)$.
\end{thm}
\begin{proof}
 It is clear from the definition that $\cals(j,k)$ has exactly one marker at every crossing point, and thus Corollary~\ref{cor 43} implies that $\cals(j,k) $ is a Kauffman state. 
 
 To show that $\cals(j,k)$ is join irreducible, suppose that it has a descent given by the inverse transposition at some segment $\ell$. Let $x$ and $y$ be the endpoints of $\ell$, and denote the level of $\ell$ by $d$. 
 Then the marker of $\cals(j,k)$ is counterclockwise from $\ell$ at $x$ and at $y$, and $\ell\in K(d)$.
 
 By definition, $K(d)$ is the closure of the set of all segments that are incident to a segment of $K(d+1)$.  The closure is constructed in a recursive manner at one crossing point at a time: if segment $a$ is in the set with endpoint $x$, we adjoin the  successors of $a$ at $x$ up to the minimal state marker. Then we repeat the construction at the other endpoints of the newly added segments, and so on. At every crossing point $x$, we can distinguish between the segments that were added to  $K(d)$ by closing at a different crossing point, and those that were added to $K(d)$ by closing at the point $x$. We call the former type \emph{incoming segments} at $x$. 
 
 Since the closing operation adds only segments that are between an incoming segment and the minimal state marker in the clockwise direction, we see that the first segment in clockwise order must be incoming. In the four cases in Figure~\ref{fig markers}, the first incoming segment is $a_1, a_2, a_3$ and $a_4$, respectively.
 In particular, the marker of the state $\cals(j,k)$ is always counterclockwise from the first incoming segment.
 
 Now let's consider the endpoints $x,y$ of our segment $\ell$ above. We have seen that the marker of $\cals(j,k)$ is counterclockwise from $\ell$ at both endpoints. Thus the segment $\ell $ must be incoming at both of its endpoints. This is impossible, unless $\ell$ is the initial segment $j$. We have shown that there is a unique segment $j$ for which the inverse transposition is a descent from $\cals(j,k)$. Thus $\cals(j,k)$ is join irreducible.

 It only remains to show that $\cals(j,1)<\cals(j,2)<\cdots <\cals(j,d_j)$. From Proposition~\ref{prop:PropertiesLattice}(c), we already know that these states are linearly ordered. 
 Consider a path $w_1$ in $L(i)$ from $\cals(j,1)$ to the minimal state. It is given by a sequence of transpositions, where the first one is at $j$ and the following ones are at the other segments in the set $K(1)$. 
 In particular $\cals(j,1) $ is the unique state that is join irreducible with descent at $j$ for which every path $w_1$ to the minimal state uses the transposition at $j$ exactly once. The set $K(2)$ is the closure of the set of segments incident to the set $K(1)$. 
 Thus there is a path $w_2$ from $\cals(j,2)$ to $\cals(j,1)$ that has an initial segment given by the same sequence of transpositions as the path $w_1$, and then uses the transpositions at the segments in $K(2)\setminus K(1)$. Similarly, there is a path $w_k$ from $\cals(j,k)$ to $\cals(j,k-1)$ that uses  the same sequence of transpositions as the path $w_{2}$. In particular, $\cals(j,k-1)<\cals(j,k)$ for all $k$.
\end{proof}



\subsection*{Acknowledgements}
 We thank Emily Gunawan for inspiring discussions on this topic.

\end{document}